\documentclass[12pt]{article}
\usepackage{latexsym, amssymb}

\newcommand{\rr}{\mathbb R}

\newcommand{\nn}{\mathbb N}

\newtheorem{theorem}{Theorem}
\newtheorem{prop}{Proposition}
\newtheorem{lemma}{Lemma}
\newtheorem{defn}{Definition}

\newtheorem{corr}{Corrolary}

\newcommand{\Cl}{\mathop{\rm Cl }\nolimits}
\newcommand{\diam}{\mathop{\rm diam}\nolimits}
\newcommand{\dist}{\mathop{\raisebox{0pt}{\rm d}}\nolimits}
\newcommand{\Int}{\mathop{\rm Int }\nolimits}
\newcommand{\myrho}{\rho_{K}}
\newcommand{\mystar}{\raisebox{-3.5pt}{*}}

\title{On imbedding of closed 2-dimensional disks into $\rr^{2}$.}
\author{Polulyakh E. A.}
\date{Institute of mathematics\\
Tereshchenkivska 3, Kiev, Ukraine.\\
e-mail: polulyah@imath.kiev.ukraine
}

\begin{document}

\maketitle

\begin{abstract}
A criterium of a closed subset of $\rr^{2}$ to be homeomorphic to a closed
2--dimensional disk using only the local information about boundary of
this subset is given in the paper.
\end{abstract}

\vspace{1ex}

It is well-known that a fiber bundle over a sircle with the Kantor set as
a fiber is called Pontryagin bundle. Studying embeddings of Pontryagin
bundles into 2--dimensional manifolds autor required  a criterium of a
closed subset of $\rr^{2}$ to be homeomorphic to a closed 2--dimensional
disk using only the local information about boundary of this subset.
The criterium of a sort is considered in what follows.

\begin{defn}
{\em
Let $X$ be a topological space, $U$ -- opened subset of $X$. We will say
that point $x \in \partial U$ is {\it accessible} from $U$ if there exists
continuous injective mapping $\varphi : I \rightarrow \Cl D$ such that
$\varphi(1)=x$, $\varphi([0,1)) \subset \Int U$.
}
\end{defn}

\begin{theorem}
{\em
Let $D$ be a compact subset of $\rr^{2}$ with nonempty interior.
$D$ is homeomorphic to a closed 2-dimentional disk if and only if the
following conditions are realized:
\begin{itemize}
        \item[1)] set $\Int D$ is connected;
        \item[2)] set $\rr^{2} \setminus D$ is connected;
        \item[3)] any $x \in \partial D$ is accessible from $\Int D$;
        \item[4)] any $x \in \partial D$ is accessible from $\rr^{2}
                  \setminus D$.
\end{itemize}
}
\end{theorem}

\noindent
{\bf Proof.} The immediate verification shows us that 2-dimentional closed
disk $D$ embedded into $\rr^{2}$ complies with the conditions 1) -- 4).

Now we are going to prove the inverse statement which appears to be
considerably more untrivial.

Let $D \subset \rr^{2}$ be a compact set satisfying the conditions 1) -- 4).

The idea of later considerations is to decompose $\partial D$ into two
parts intersecting at two points and to verify the gomeomorphism of each
of the parts to an interval. That will imply the gomeomorphism of
$\partial D$ to a sircle $S^{1}$ which involves the gomeomorphism of $D$
to a disk.

First of all we will prove two simple statements we will make use of to
determinate main construction.

\begin{lemma}
{\em
Let $U$ be an opened connected subset of $\rr^{2}$. For any $x_{1}$,
$x_{2} \in U$ there exists a continuous injective mapping $\psi:I
\rightarrow U$ such that $\psi(0)=x_{1}$, $\psi(1)=x_{2}$.
\label{lemma1}
}
\end{lemma}

\noindent
{\bf Proof.} It is known that an opened connected subset of locally arcwise
connected space is arcwise connected. This implies the existance of
continuous mapping
$$\psi':I \rightarrow U ,$$
such that $\psi'(0)=x_{1}$, $\psi'(1)=x_{2}$.

For each $y \in \psi'(I)$ we find $\varepsilon = \varepsilon(y) > 0$ to
satisfy the implication $B_{\varepsilon}(y) \subset U$. The collection
$\{ B_{\varepsilon}(y) \}_{y \in \psi'(I)}$ form the opened covering of
$\psi'(I)$. Since the set $\psi'(I)$ is compact we can select a finit
subcovering $\{ W_{i} \}_{i=1}^{n}$.

For $i,j \in \{ 1, \ldots, n \}$, $i \ne j$ we join centers of disks
$W_{i}$ and $W_{j}$ having nonempty intersection by closed interval (it is
evident that the interval belongs to $W_{i} \cup W_{j}$). We add
intervals which connect points $x_{1}$, $x_{2}$ with centers of disks from
the collection $\{ W_{i} \}_{i=1}^{n}$ containing these points to the
system of intervals received.

We obtain intervals from the system we built to be mutually nonparallel by
small perturbation of centers and radii of disks $W_{i}$, $i=1, \ldots, n$
remaining the properties of $\{ W_{i} \}_{i=1}^{n}$ to be a covering of
$\psi'(I)$ and to belong to $U$.

After we add all intersections of intervals to the set containing points
$x_{1}$, $x_{2}$ and the senters of disks $W_{i}$, $i=1, \ldots, n$ we get
a finit point set $l$. Denote by $L$ the system of closed intervals with
endpoints in $l$ which complies with the following conditions:
\begin{enumerate}
        \item[---] every element of $L$ is subset of certain interval from
                   the system of intervals built beyond;
        \item[---] the intersection of interior with the set $l$ is empty
                   for each element of $L$.
\end{enumerate}
We receive finite graph
$$R=(l, L)$$
embedded into $\rr^{2}$ with the set of vertices $l$ and the set of edges
$L$. Moreover we have an implication
$$R \subset \bigcup_{i=1}^{n} W_{i} \;.$$

The graph $R$ is connected. Differently the set $\bigcup_{i=1}^{n} W_{i}$
have to decompose into two nonintersecting opened sets $W^{(1)}$,
$W^{(2)}$ realizing the equations
$$\psi'(I) \subset W^{(1)} \cup W^{(2)} \;,\;
\psi'(I) \cap W^{(1)} \ne \emptyset \;,\;
\psi'(I) \cap W^{(2)} \ne \emptyset \;.$$
contradicting to arcwise connectivity of $\psi'(I)$.

It is known that connected graph $R$ can be reduced by removing of certain
edges to a tree $\widetilde{R} = (l, \widetilde{L})$, $\widetilde{L}
\subset L$. It is known as well that every two vertices of a tree could be
connected by a simple path.

Let us take the simple path connecting vertices $x_{1}$ and $x_{2}$ of
$\widetilde{R}$ and construct continuous one-to-one mapping $\psi$ of the
interval $I$ onto this path. We have continuous injective mapping
$$\psi:I \rightarrow U \;,\; \psi(0) = x_{1} \;,\; \psi(1) = x_{2} \;.$$
Q. E. D.

\begin{lemma}
{\em
Let $U$ be an opened connected subset of $\rr^{2}$; $x_{1}$, $x_{2} \in
\partial U$ be accessible from $U$, $x_{1} \ne x_{2}$. There exists a
continuous injective mapping $\psi:I \rightarrow \Cl U$ such that
$\psi(0)=x_{1}$, $\psi(1)=x_{2}$, $\psi((0, 1)) \subset U$.
\label{lemma2}
}
\end{lemma}

\noindent
{\bf Proof.} Let us take continuous injective mappings $\psi_{s}:I
\rightarrow \Cl U$, $s=1, 2$ realizing the conditions $\psi_{s}(0)=x_{s}$,
$\psi_{s}((0,1]) \subset U$, $s=1, 2$.

Mappings $\psi_{1}$, $\psi_{2}$ could be selected so as to satisfy
equality $\psi_{1}(I) \cap \psi_{2}(I) = \emptyset$. Really, assume it is
not the case. We denote
$$t_{0} = \min \{ t \in I \; | \; \psi_{1}(t) \in \psi_{2}(I) \} .$$
Since the compact sets $x_{1}=\psi_{1}(0)$ and $\psi_{2}(I)$ have a trivial
intersection (i. e. $\dist (x_{1}, \psi_{2}(I)) > 0$) and to addition
mapping $\psi_{1}$ is continuous, the inequality $t_{0} > 0$ is valid and
$\psi_{1}$ could be replaced by mapping
$$\widetilde{\psi}_{1}:I \rightarrow \Cl U ,$$
$$\widetilde{\psi}_{1}(t) = \psi_{1}(\frac{t t_{0}}{2}),$$
so that $\widetilde{\psi}_{1}(I) \cap \psi_{2}(I) = \emptyset$.

Denote $z_{1}=\psi_{1}(1) \in U$, $z_{2}=\psi_{2}(1) \in U$. In accord
with lemma \ref{lemma1} there exists a continuous injective mapping $\psi':I
\rightarrow U$ such that $\psi'(0)=z_{1}$, $\psi'(1)=z_{2}$.

Assume
$$t_{s} = \min \{ t \in I \;|\; \psi_{s}(t) \in \psi'(I) \}
\;,\; s=1,2 \;.$$
Since the compact sets $x_{s}$ and $\psi'(I)$ are mutually disjoint (i. e.
$\dist (x_{s}, \psi'(I)) > 0$), $s=1, 2$ and to addition the mappings
$\psi_{s}$ are continuous, the inequalities $t_{s} > 0$, $s=1, 2$ are
realized. As mapping $\psi'$ is injective, points $\tau_{s} \in I$, $s=1,
2$ complying with equalities $\psi'(\tau_{s}) = \psi_{s}(t_{s})$ are
unumbiguously defined.

Under $\tau_{1} < \tau_{2}$ we denote
$$\alpha=t_{1} + t_{2} + \tau_{2} - \tau_{1},$$
$$\psi(t)=\cases{
                    \psi_{1}(\alpha t) & when $t \in [0,
                          \frac{t_{1}}{\alpha})$ \cr
                    \psi'(\alpha t - t_{1} + \tau_{1}) & when
                          $t \in [\frac{t_{1}}{\alpha},
                          \frac{t_{1}+\tau_{2}-\tau_{1}}{\alpha}]$ \cr
                    \psi_{2}(\alpha-\alpha t) & when
                          $t \in (\frac{t_{1}+\tau_{2}-
                          \tau_{1}}{\alpha},1]$ \cr
                 } \;,$$
else
$$\alpha=t_{1} + t_{2} + \tau_{1} - \tau_{2},$$
$$\psi(t)=\cases{
                    \psi_{1}(\alpha t) & when $t \in [0,
                          \frac{t_{1}}{\alpha})$ \cr
                    \psi'(\alpha t - t_{1} + \tau_{2}) & when
                          $t \in [\frac{t_{1}}{\alpha},
                          \frac{t_{1}+\tau_{1}-\tau_{2}}{\alpha}]$ \cr
                    \psi_{2}(\alpha-\alpha t) & when
                          $t \in (\frac{t_{1}+\tau_{1}-
                          \tau_{2}}{\alpha},1]$ \cr
                 } \;.$$
Mapping $\psi:I \rightarrow I^{2}$ meets with the conclusion of lemma
\ref{lemma2}. \\
Lemma is proved.

\bigskip

The set $\partial D$ divides $\rr^{2}$ into two connected components, thus
the dimension of $\partial D$ is not less than one [1] and it containes more
than one point.

We fix two distinct points $z_{1}$, $z_{2} \in \partial D$. According to
lemma \ref{lemma2} there exist continuous injective mappings
$$\varphi_{1}:I \rightarrow D, \; \varphi_{2}:I \rightarrow \rr^{2}
\setminus \Int D$$
such that $\varphi_{1}((0,1)) \subset \Int D$, $\varphi_{2}((0,1)) \subset
\rr^{2} \setminus D$, $\varphi_{i}(0) = z_{1}$, $\varphi_{i}(1) = z_{2}$
$i=1,2$.

The mapping
$$\gamma:I \rightarrow \rr^{2}$$
$$\gamma(t) = \cases{
                    \varphi_{1}(2 t) & when $t \in [0,
                          \frac{1}{2})$ \cr
                    \varphi_{2}(2(1-t)) & when
                          $t \in [\frac{1}{2},1]$ \cr
                 }$$
could be considered as the mapping $\gamma:S^{1}=I/(\{ 0 \} \cup \{ 1 \})
\rightarrow \rr^{2}$ being continuous and injective.

It is known that a sircle imbedded into $\rr^{2}$ bounds a disk.
Consequently the set $\gamma (S^{1})$ bounds a closed subset $B_{1}$ of
$\rr^{2}$ homeomorphic to disk.

We have already proved that $\gamma (S^{1})$ decomposes $\rr^{2}$ into two
connected components $\Int B_{1}$ and $\rr^{2} \setminus B_{1}$, with any
point of $\gamma (S^{1}) = \partial B_{1}$ being accessible both from $\Int
B_{1}$ and from $\rr^{2} \setminus B_{1}$.

Let us show that $\partial D \cap \Int B_{1} \ne \emptyset$. To this end
we fix $y_{1} \in \varphi_{1}((0,1)) \subset \Int D$, $y_{2} \in
\varphi_{2}((0,1)) \subset \rr^{2} \setminus D$. Lemma \ref{lemma2}
guarantees the existence of a continuous injective mapping $\psi:I
\rightarrow B_{1}$ such that $\psi(0)=y_{1}$, $\psi(1)=y_{2}$,
$\psi((0,1)) \subset \Int B_{1}$. Since the points $y_{1}$, $y_{2}$
belongs to distinct connected components of $\rr^{2} \setminus \partial D$
then $\psi(I) \cap \partial D \ne \emptyset$ and there exists $t \in (0,1)$
such that $\psi(t) \in \partial D$.

Similarly is proved that $\partial D \cap ( \rr^{2} \setminus B_{1} ) \ne
\emptyset$.

\begin{lemma}
{\em
The opened sets $\widetilde{D}^{(1)}_{1} = \Int B_{1} \cap \Int D$,
$\widetilde{D}^{(1)}_{2} = \Int B_{1} \cap (\rr^{2} \setminus D)$ are
arcwise connected.
}\label{lemma3}
\end{lemma}

{\bf Proof.} Let $y_{1}$, $y_{2} \in \widetilde{D}_{1}^{(1)}$. In accord
with lemma \ref{lemma1} there exists an injective continuous mapping
$\psi:I \rightarrow \Int D$, $\psi(0)=y_{1}$, $\psi(1)=y_{2}$. If $\psi(I)
\cap \partial B_{1} = \emptyset$ then $\psi(I) \subset \Int B_{1}$ and the
points $y_{1}$, $y_{2}$ are contained in the same connected component of
$\widetilde{D}^{(1)}_{1}$.

Suppose $\psi(I) \cap \partial B_{1} \ne \emptyset$. We assume
$$t_{1} = \min \{ t \in I \; | \; \psi(t) \in \partial B_{1} \} \;,$$
$$t_{2} = \max \{ t \in I \; | \; \psi(t) \in \partial B_{1} \} \;.$$
Since $y_{1} = \psi(0) \in \Int B_{1}$ and $y_{2} = \psi(1) \in \Int
B_{1}$, the values $t_{1}$, $t_{2}$ meet an inequality $0 < t_{1} \le
t_{2} < 1$.

Consider the unumbigously defined values $\tau_{1}$, $\tau_{2} \in (0,1)$
complying with the equalities $\varphi_{1}(\tau_{i}) = \psi(t_{i})$,
$i=1,2$. For defenitness we can suppose $\tau_{1} \le \tau_{2}$. Compacts
$\partial D$ and $\varphi_{1}([\tau_{1}, \tau_{2}])$ have a trivial
intersection, thus
$$\varepsilon = \dist (\partial D, \varphi_{1}([\tau_{1}, \tau_{2}])) >
0.$$

Let us fix a homeomorphism
$$f:B_{1} \rightarrow \{ (r, \varphi) \in \rr^{2} \; | \; | r | \le 1
\}.$$
Mappings $f$, $f^{-1}$ are uniformly continuous. Therefore $\delta > 0$
could be found to realize an implication of $\dist (f^{-1}(z_{1}),
f^{-1}(z_{2})) < \varepsilon / 2$ from $\dist (z_{1}, z_{2}) < \delta$.
Also $\delta_{0} > 0$ could be found such that $\dist (x_{1}, x_{2}) <
\delta_{0}$ implies $\dist (f(x_{1}), f(x_{2})) < \delta$.

Since the mapping $\psi:I \rightarrow \Int B_{1}$ is continuous, there
exists $t_{1}'$, $t_{2}' \in (0,1)$, $0 < t_{1}' < t_{1}$, $t_{2} < t_{2}'
< 1$ satisfying to the inequalities
$$\dist (\psi(t_{i}), \psi(t_{i}')) < \delta_{0}, \; i=1,2 .$$

Under the mapping $f$ the set $\varphi_{1}([\tau_{1}, \tau_{2}])$ pases to
an arc
$$S_{0} = \{ (r, \varphi) \in \rr^{2} \; | \; r=1, \; \varphi \in
[\varphi_{1}, \varphi_{2}] \}$$
of the sircle $\{ (r, \varphi) \in \rr^{2} \; | \; r=1 \}$. Here
$\varphi_{i}$ is a certain polar angle corresponding to the point
$f \circ \varphi_{1}(\tau_{i})$, $i=1,2$.

There exists a continuous mapping $\widetilde{\alpha}:I \rightarrow \{ (r,
\varphi) \in \rr^{2} \; | \; |r| < 1 \}$, $\widetilde{\alpha}(0) = f \circ
\psi(t_{1}')$, $\widetilde{\alpha}(1) = f \circ \psi(t_{2}')$ such that
$\widetilde{\alpha}(I) \subset U_{\delta}(S_{0})$. We can assume for
instance
$$\widetilde{\alpha}(t) = (tr_{2}' + (1-t)r_{1}', t \varphi_{2}' + (1-t)
\varphi_{1}') ,$$
where $(r_{i}', \varphi_{i}')$ is an appropriate presentation of $f \circ
\psi (t_{i}')$. $i=1,2$ in polar coordinats.

An image of the mapping $\alpha = f^{-1} \circ \widetilde{\alpha}:I
\rightarrow \Int B_{1}$ lies in the $\varepsilon / 2$ -- neighbourhood of
the arc $\varphi_{1}([\tau_{1}, \tau_{2}])$, thus $\alpha(I) \cap \partial
D = \emptyset$.

Denote
$$a = 2 + t_{1}' - t_{2}' ,$$
$$\widetilde{\psi}(t) = \cases{
                    \psi(at) & when $t \in [0, \frac{t_{1}'}{a})$ \cr
                    \alpha(at-t_{1}') & when $t \in [\frac{t_{1}'}{a},
                          \frac{1+t_{1}'}{a})$ \cr
                    \psi(at-1-t_{1}'+t_{2}') & when
                          $t \in [\frac{1+t_{1}'}{a},1]$ \cr
                 }$$
We obtain the continuous mapping $\widetilde{\psi}:I \rightarrow \Int B_{1}$
such that $\widetilde{\psi}(0) = y_{1}$, $\widetilde{\psi}(1) =
y_{2}$, $\widetilde{\psi}(I) = \psi([0, t_{1}') \cup (t_{2}', 1]) \cup
\alpha(I)$. Therefore points $y_{1}$, $y_{2}$ belong to a same connected
component of $\widetilde{D}^{(1)}_{1}$. The arcwise connectedness of
$\widetilde{D}^{(1)}_{1}$ follows from the arbitrary rool of selection for
$y_{1}$, $y_{2}$.

Analogously is proved the arcwise connectedness of
$\widetilde{D}^{(1)}_{2}$. \\
Lemma \ref{lemma3} is proved.

\bigskip

We arive to the following construction. Let
$$K^{(1)} = \partial D \cap B_{1}, \; D^{(1)}_{1} = B_{1} \cap \Int D, \;
D^{(1)}_{2} = B_{1} \cap (\rr^{2} \setminus D).$$
At that time $B_{1} = D^{(1)}_{1} \cup D^{(1)}_{2} \cup K^{(1)}$ where
$D^{(1)}_{1}$, $D^{(1)}_{2}$ are opened in $B_{1}$ arcwise connected sets
complying with the equalities $D^{(1)}_{1} \cap D^{(1)}_{2} = \emptyset$,
$K^{(1)} = \partial D^{(1)}_{1} = \partial D^{(1)}_{2}$. Any point of
$K^{(1)}$ is accesible both from $D^{(1)}_{1}$ and $D^{(1)}_{2}$.
Furthermore $K^{(1)} \cap \partial B_{1} = \{ z_{1} \} \cup \{ z_{2} \}$
for certain $z_{1} \ne z_{2}$.

Let us construct an analogous design containing a homeomorphic image of
the set $\partial D \cap (\rr^{2} \setminus \Int B_{1})$.

We fix a point $z_{0} \in \Int B_{1} \cap \Int D_{1}$. There exists a
linear automorphism of $\rr^{2}$ translating $z_{0}$ to the origin of
coordinates. Thus we can assume the origin to be contained in $\Int B_{1}
\cap \Int D_{1}$.

Consider the involution
$$f: \rr^{2} \setminus \{ 0 \} \rightarrow \rr^{2} \setminus \{ 0 \}, \;
f(r, \varphi) = (r^{-1}, \varphi),$$
which is known to be an automorphism of $\rr^{2} \setminus \{ 0 \}$.

Note the opened set $(\Int D) \setminus \{ 0 \}$ is connected. Really, fix
$y_{1}$, $y_{2} \in \Int D \setminus \{ 0 \}$. We shall take $\psi:I
\rightarrow \Int D$ -- a continuous curve in $\Int D$ connecting $y_{1}$
and $y_{2}$. Since $0 \in \Int D$ there exists $\delta > 0$ such that
$U_{\delta}(0) \subset \Int D$. Let
$$\varepsilon = \frac{1}{2} \min (\delta, \dist(y_{1}, 0), \dist(y_{2},
0)).$$
It is evident that $\Cl U_{\varepsilon}(0) \subset \Int D$. In the case
$\psi(I) \cap \Cl U_{\varepsilon}(0) = \emptyset$ points $y_{1}$, $y_{2}$
are contained in a same connected component of $\Int D \setminus \{ 0 \}$.
Presuppose $\psi(I) \cap \Cl U_{\varepsilon}(0) \ne \emptyset$.

Denote
$$t_{1} = \min \{ t \in I \; | \; \dist (\psi(t), 0) \le \varepsilon \},$$
$$t_{2} = \max \{ t \in I \; | \; \dist (\psi(t), 0) \le \varepsilon \}.$$
By virtue of choise of $\varepsilon$ the inequalities $0 < t_{1}, t_{2} <
1$ are valid. Let us consider an arc of the sircle bounding
$U_{\varepsilon}(0)$ with endpoints $\psi(t_{1})$, $\psi(t_{2})$.
Substituting the curve $\psi: [t_{1}, t_{2}] \rightarrow \Int D$ with this
arc we receive a continuous path $\widetilde{\psi}:I \rightarrow \Int D
\setminus \{ 0 \}$ joining $y_{1}$ with $y_{2}$. On the strength of
arbitrary rool for selection of $y_{1}$, $y_{2}$ the set $\Int D \setminus
\{ 0 \}$ is connected.

Hence the set $\rr^{2} \setminus (\{ 0 \} \cup f(\partial D))$ have two
connected components.

Since the set $D$ is limited and
relation $\stackrel{\circ}{U}_{\varepsilon}(0) = \{ x \in \rr^{2} \; | \;
\dist (x, 0) \in (0, \varepsilon) \} = f(\{ x \in \rr^{2} \; | \; \dist
(x, 0) > 1/ \varepsilon \})$ is valid for all $\varepsilon > 0$, there
exists $\delta > 0$ such that $\stackrel{\circ}{U}_{\delta}(0) \subset
f(\rr^{2} \setminus D) \subset \rr^{2} \setminus f(\partial D)$. Therefore
in the first place $0 \in \Int (\rr^{2} \setminus f(\partial D))$ and the
set $\rr^{2} \setminus f(\partial D)$ is opened; in the second $0$ and
$f(\rr^{2} \setminus D)$ are contained in the same connected component of
$\rr^{2} \setminus f(\partial D)$ while $0$ and $f(\Int D)$ -- in distinct
components. Consequently the set $\rr^{2} \setminus f(\partial D)$ have
two connected components one of which is limited since $U_{\varepsilon}(0)
\subset \Int D$ for certain $\varepsilon > 0$ and $f(U_{\varepsilon}(0)) =
\{ x \in \rr^{2}\; | \; \dist(x, 0) \ge 1 / \varepsilon \} \subset f(\Int
D)$. Denote this component by $V$. At that time $f(\rr^{2} \setminus D)
\subset V$ and $\partial V = f(\partial D)$.

Denote $D' = \Cl V$. we received the compact subset of $\rr^{2}$
complying with conditions 1) -- 4) of the theorem with any point of
$\partial D$ being accessible even from $f(\rr^{2} \setminus D) \subset
\Int D'$ and $f(\Int D) = \rr^{2} \setminus D'$).

Consider continuous injective mappings
$$\widetilde{\varphi}_{1} =
f \circ \varphi_{1}:I \rightarrow \rr^{2} \setminus \Int D', \;
\widetilde{\varphi}_{2} = f \circ \varphi_{2}:I \rightarrow D'.$$
They meet with the conditions $\widetilde{\varphi}_{1}((0, 1)) \subset
\rr^{2} \setminus D'$,
$\widetilde{\varphi}_{2}((0, 1)) \subset \Int D'$,
$\widetilde{\varphi}_{i}(0) = f(z_{1}) \in \partial D'$,
$\widetilde{\varphi}_{i}(1) = f(z_{2}) \subset \partial D'$, $i=1,2$.

The mapping
$$\widetilde{\gamma}:S^{1} = I / (\{ 0 \} \cup \{ 1 \}) \rightarrow
\rr^{2},$$
$$\widetilde{\gamma}(t) = \cases{
                    \widetilde{\varphi_{1}}(2t) & when $t \in
                          [0, \frac{1}{2})$ \cr
                    \widetilde{\varphi_{2}}(2(1-t)) & when $t \in
                          [\frac{1}{2}, 1)$ \cr
                 }$$
appears to be continuous and injective. Therefore the set
$\widetilde{\gamma}(S^{1})$ bounds a closed disk $B_{2} \subset \rr^{2}$.

The immediate verification shows that
$f(\rr^{2} \setminus \Int B_{1}) \subset B_{2}$, $f(B_{1}) = \rr^{2}
\setminus \Int B_{2}$, thus $K^{(2)} = f(\partial D \setminus \Int B_{1}) =
\partial D' \cap B_{2}$.

Denote
$$D^{(2)}_{1} = B_{2} \cap \Int D', \; D^{(2)}_{2} = B_{2} \cap (\rr^{2}
\setminus D').$$
By repeating argument mentioned beyond for $D^{(1)}_{1}$, $D^{(1)}_{2}$ we
draw a conclusion that $D^{(2)}_{1}$, $D^{(2)}_{2}$ are opened in $B_{2}$
arcwise connected sets satisfying the equations $B_{2} = D^{(2)}_{1} \cup
D^{(2)}_{2} \cup K^{(2)}$, $D^{(2)}_{1} \cap D^{(2)}_{2} = \emptyset$,
$K^{(2)} = \partial D^{(2)}_{1} = \partial D^{(2)}_{2}$, with any point of
$K^{(2)}$ being accessible both from $D^{(2)}_{1}$ and $D^{(2)}_{2}$.
Furthemore $K^{(2)} \cap \partial B_{2} = \{ f(z_{1}) \} \cup \{ f(z_{2}) \}$
($f(z_{1}) \ne f(z_{2})$.

Consider the space $(I^{2},\dist)$, where $I^{2}=[0,1] \times [0,1]$ and
$\dist$ is a metric on $I^{2}$.

Consider also a triple $D_{1}$, $D_{2}$, $K$ complying with the following
conditions. $I^{2}=D_{1} \cup D_{2} \cup K$, where $D_{1}$, $D_{2}$ are
opened in $I^{2}$ arcwise connected sets realizing the equalities
$D_{1} \cap D_{2} = \emptyset$, $K = \partial D_{1} = \partial D_{2}$,
with any point of $K$ being accessible both from $D_{1}$ and $D_{2}$.
Furthemore $K \cap \partial I^{2} = \{ a \} \cup \{ b \}$ for certain $a
\in \{ 0 \} \times (0,1)$, $b \in \{ 1 \} \times (0,1)$.

There exist homeomorphisms $h_{i}:B_{i} \rightarrow I^{2}$, $i=1,2$ such
that $h_{1}(z_{1})$, $h_{2} \circ f(z_{1}) \in \{ 0 \} \times [0,1]$;
$h_{1}(z_{2})$, $h_{2} \circ f(z_{2}) \in \{ 1 \} \times [0,1]$.

Under the mapping $h_{i}$ the triple $D^{(i)}_{1}$, $D^{(i)}_{2}$,
$K^{(i)}$ passes to a triple kind $D_{1}$, $D_{2}$, $K$.

To complete the proof of theorem it is sufficient to verify the set $K$
from the triple $D_{1}$, $D_{2}$, $K$ specified above is homeomorphic to
an interval.

Denote
$$ \{ 0 \} \times [0,1] = I_{l} \;,\; \{ 0 \} \times (0,1) = \stackrel{\circ}
{I_{l}} \;,$$
$$ \{ 1 \} \times [0,1] = I_{r} \;,\; \{ 1 \} \times (0,1) = \stackrel{\circ}
{I_{r}} \;,$$
$$ [0,1] \times \{ 0 \} = I_{b} \;,\; (0,1) \times \{ 0 \} = \stackrel{\circ}
{I_{b}} \;,$$
$$ [0,1] \times \{ 1 \} = I_{t} \;,\; (0,1) \times \{ 1 \} = \stackrel{\circ}
{I_{t}} \;.$$

\begin{prop}
{\em
Let $x$ be contained in $K \setminus (\{ a \} \cup \{ b \})$. There exists
a continuous injective mapping $\varphi_{x} : I \rightarrow I^{2}$
complying with relations
$$ \varphi_{x}(0) \in \stackrel{\circ}{I_{b}} \;,\; \varphi_{x}(1) \in
\stackrel{\circ}{I_{t}} \;,\; \varphi_{x}(1/2) = x \;,$$
$$ \varphi_{x}(0,1) \subset  \Int I^{2} \;,\;  \varphi_{x}([0, 1/2))
\subset  D_{1}   \;,\;  \varphi_{x}((1/2,1])   \subset  D_{2}   \;.$$
\label{prop1}
}
\end{prop}

\noindent
{\bf Proof.} Examine the opened sets
$$U_{1} = D_{1} \cap (I^{2} \setminus \partial I^{2}), \;
U_{2} = D_{2} \cup (I^{2} \setminus \partial I^{2}).$$
By analogy with lemma \ref{lemma3} it is proved that the sets $U_{1}$, $U_{2}$
are arcwise connected.

We fix the points $z_{1} \in \stackrel{\circ} {I_{b}}$, $z_{2} \in
\stackrel{\circ} {I_{t}}$. Any of these points is accessible from $I^{2}
\setminus \partial I^{2}$ since $I^{2}$ is homeomorphic to a disk.
Consequently continuous injective mappings
$$\varphi_{i}':I \rightarrow I^{2}, \; i=1,2$$
could be found to realize the correlations $\varphi_{i}'(0)=z_{i}$,
$\varphi_{i}'((0,1]) \in I^{2} \setminus \partial I^{2}$.

Denote
$$t_{i} = \min \{ t \in I \; | \; \varphi_{i}'(t) \in K \}, \; i=1,2 .$$
As $z_{i} = \varphi_{i}'(0) \not\in K$ then $t_{i} > 0$, $i=1,2$.

It is obvious that injective continuous mappings
$$\varphi_{i}:I \rightarrow I^{2}, \; \varphi_{i}(t) =
\varphi_{i}'(\frac{t t_{i}}{2}), \; i=1,2$$
meet the conditions $\varphi_{i}(0)=z_{i}$, $\varphi_{i}((0,1]) =
\varphi_{i}'((0, t_{i}/2]) \subset (I^{2} \setminus
\partial I^{2}) \cap D_{i} = U_{i}$, $i=1,2$.

Hence points $z_{i}$ are accessible from $U_{i}$, $i=1,2$.

Similarly is proved that $x \in K \setminus (\{ a \} \cup \{ b \})$ is
accessible both from $U_{1}$ and $U_{2}$.

In accord with lemma \ref{lemma2} there exist continuous injective mappings
$$\psi_{i}:I \rightarrow \Cl U_{i} = \Cl D_{i}, \; i=1,2$$
such that $\psi_{i}(0)=z_{i}$, $\psi_{i}(1)=x$, $\psi_{i}((0,1)) \subset
U_{i} \subset D_{i}$.

Since $D_{1} \cap D_{2} = K = \partial D_{1} = \partial D_{2}$, $U_{i}
\subset \Int D_{i}$, $i=1,2$, then $U_{1} \cap U_{2} = \emptyset$ and the
continuous mapping
$$\varphi_{x}:I \rightarrow I^{2},$$
$$\varphi_{x}(t) = \cases{
                   \psi_{1}(2t), & for $t \in [0,\frac{1}{2})$ \cr
                   \psi_{2}(2(1-t)), & for $t \in
                   [\frac{1}{2}, 1]$ \cr
}$$
appeares to be injective and it complies with all conditions of the
proposition.

\begin{prop}
{\em
Let $x,y$, $K \setminus (\{ a \} \cup \{ b \})$, $x \ne y$, $\varphi_{x}:I
\rightarrow I^{2}$ be a mapping satisfying the conditions of
proposition \ref{prop1}. There exists a continuous injective mapping
$\varphi_{y}:I \rightarrow I^{2}$ realizing conditions of proposition
\ref{prop1} and such that $\varphi_{x}(I) \cap \varphi_{y}(I) =
\emptyset$.
\label{prop2}
}
\end{prop}

\noindent
{\bf Proof.} For defenitness we suppose that $y$ and $I_{r}$ are contained
in the same connected component of  $I^{2} \setminus \varphi_{x}(I)$.

Consider the disk $B$ bounded by $\varphi_{x}(I) \cup ([\varphi_{x}(1), 1]
\times \{ 1 \}) \cup I_{r} \cup ([\varphi_{x}(0), 1] \times \{ 0 \})$.
There exists a homeomorphism $f:B \rightarrow I^{2}$ meeting the
equalities $f(\varphi_{x}(I)) = I_{l}$, $f(I_{r}) = I_{r}$.

Remark that $V_{1} = f(D_{1} \cap B)$ and $V_{2} = f(D_{2} \cap B)$ are
opened in $I^{2}$ nonintersecting arcwise connected sets with the common
boundary $\widetilde{K} = f(K \cap B)$ every point of which is accessible
both from $V_{1}$ and $V_{2}$.

According to proposition \ref{prop1} a continuous injective mapping
$\psi_{f(y)}: I \rightarrow I^{2}$ could be found to comply with relations
$$
\begin{array}{c}
\psi_{f(y)}(0) \in \stackrel{\circ} {I_{b}}\,, \quad
\psi_{f(y)}(1) \in \stackrel{\circ}{I_{t}}\,, \quad
\psi_{f(y)}(1/2) = f(y)\,, \vphantom{\sum\limits_{1}^{1}}\\
\psi_{f(y)}((0,1)) \subset \Int I^{2}\,, \quad
\psi_{f(y)}[0, 1/2) \subset V_{1}\,, \quad
\psi_{f(y)}(1/2, 1] \subset V_{2}\,.
\end{array}
$$

The mapping $\varphi_{y} = f^{-1} \circ \psi_{f(y)}:I \rightarrow I^{2}$
will be desired.

\bigskip

We will define a linear order relation on $K$.

\begin{defn}
{\em
Let $x$, $y$ be contained in $K \setminus (\{ a \} \cup \{ b \})$,
$\varphi_{x}:I \rightarrow I^{2}$ be a mapping satisfying the conditions
of proposition ref{prop1}. It divides $I^{2}$ into two closed disks $H_{1}$ and
$H_{2}$, so as $H_{1} \cap H_{2} = \varphi_{x}(I)$, $I_{l} \subset H_{1}$,
$I_{r} \subset H_{2}$. We will say $y \le x$ if $y \in H_{1}$, $y \ge x$
if $y \in H_{2}$. We set $a \le x$ and $b \ge x$ for any $x \in K$.
}
\end{defn}

The following statement will be necessary in order to check the
correctness of the definition.

\begin{lemma}
{\em
Let $x,y \in \Int I^{2}$, $x \ne y$. Let $\varphi_{x}'$, $\varphi_{x}''$,
$\varphi_{y}:I \rightarrow I^{2}$ be continuous injective mappings
complying with the following conditions:

\begin{itemize}
      \item[1)] $\varphi_{x}'(0), \varphi_{x}''(0), \varphi_{y}(0) \in
                \stackrel{\circ}{I_{b}}$;
      \item[2)] $\varphi_{x}'(1), \varphi_{x}''(1), \varphi_{y}(1) \in
                \stackrel{\circ}{I_{t}}$;
      \item[3)] $\varphi_{x}'(1/2)=\varphi_{x}''(1/2)=x$, $\varphi_{y}(1/2)
                =y$;
      \item[4)] $\varphi_{x}'(0,1), \varphi_{x}''(0,1), \varphi_{y}(0,1)
                \subset I^{2} \setminus \partial I^{2}$;
      \item[5)] $\varphi_{x}'(I) \cap \varphi_{y}(I) = \emptyset$;
      \item[6)] Point $y$ and the set $I_{r}$ are contained in the same
                connected component of $I^{2} \setminus \varphi_{x}'(I)$;
      \item[7)] Point $y$ and the set $I_{l}$ are contained in the same
                connected component of $I^{2} \setminus \varphi_{x}''(I)$.
\end{itemize}
There exists a continuous mapping $\psi:I \rightarrow I^{2}$ realizing
relations
\begin{itemize}
      \item[1)] $\psi(I) \subset \varphi_{x}''(I) \cup \varphi_{y}(I)$;
      \item[2)] $\psi(0) \in \stackrel{\circ}{I_{b}}$, $\psi(1) \in
                \stackrel{\circ}{I_{t}}$;
      \item[3)] $x, y \not\in \psi(I)$.
\end{itemize}
\label{lemma4}
}
\end{lemma}

\noindent
{\bf Proof.} The set $\varphi_{x}'$ decomposes $I^{2}$ into two closed
disks $A_{1}'$ and $A_{2}'$ so that $A_{1}' \cup A_{2}' = I^{2}$,
$A_{1}' \cap A_{2}' = \varphi_{x}'(I)$, $A_{1}' \supset I_{l}$, $A_{2}'
\supset I_{r}$.

Analogously the set $\varphi_{x}''$ divides $I^{2}$ into two closed disks
$A_{1}''$ and $A_{2}''$, $A_{1}'' \cup A_{2}'' = I^{2}$,
$A_{1}'' \cap A_{2}'' = \varphi_{x}''(I)$, $A_{1}'' \supset I_{l}$, $A_{2}''
\supset I_{r}$.

The set $\varphi_{y}$ partitions $I^{2}$ into two closed disks $B_{1}$ and
$B_{2}$, $B_{1} \supset I_{l}$, $B_{2} \supset I_{r}$, with $y \in \Int
A_{1}'' \cap \Int A_{2}'$, $B_{2} \subset A_{2}'$ and $x \not\in B_{2}$.
Therefore $\varphi_{x}''((t_{1}, t_{2})) \subset B_{2}$ involves either
$t_{1}, t_{2} < 1/2$ or $t_{1}, t_{2} > 1/2$.

Let us fix a point $z \in I_{r}$. Mark $z \in \Int B_{2} \cap \Int
A_{2}''$ in the topology induced from $I^{2}$. Join the points $y$ and $z$
with an injective continuous curve $\widetilde{\gamma}:I \rightarrow B_{2}$
such that $\widetilde{\gamma}(0)=y$, $\widetilde{\gamma}(1)=z$,
$\widetilde{\gamma}((0,1)) \subset \Int B_{2}$.

Since $y \in \Int A_{1}''$ we can find $\varepsilon > 0$ to satisfy the
inclusion $B_{\varepsilon}(y) \subset \Int A_{1}''$. By virtue of the
continuity of $\widetilde{\gamma}$ there exists $t_{1} \in (0, 1)$
complying with the inclusion $\widetilde{\gamma}([0, t_{1}]) \subset
B_{\varepsilon}(y)$. Consequently the equalities $\widetilde{\gamma}([0,
t_{1}]) \cap A_{2}'' = \emptyset$, $\widetilde{\gamma}([0, t_{1}]) \cap
\varphi_{x}''(I) = \emptyset$ are valid. Analogously $t_{2} \in (0, 1)$
could be found to meet the formula $\widetilde{\gamma}([t_{2}, 1]) \cap
\varphi_{x}''(I) = \emptyset$.

Denote the following mapping by $\widetilde{\mu}:I \rightarrow B_{2}$:
$\widetilde{\mu}(t) = \widetilde{\gamma}((t_{2}-t_{1})t + t_{1})$,
$\widetilde{\mu}(0)=\widetilde{\gamma}(t_{1})$,
$\widetilde{\mu}(1)= \widetilde{\gamma}(t_{2})$,
$\widetilde{\mu}(I) \subset \Int B_{2}$.

There exists a small perturbation of $\widetilde{\mu}$ resulting in the
intersection of $\widetilde{\mu}$ and $\varphi_{x}''$ to become
transversal. Therefore the mapping $\mu:I \rightarrow B_{2}$ could be
found to meet the conditions $\mu(0)=\widetilde{\gamma}(t_{1})$,
$\mu(1)= \widetilde{\gamma}(t_{2})$, $\mu(I) \subset \Int B_{2}$ with
$\mu$ and $\varphi_{x}''$ to be intersected transversally in the finit
number of points, as the sets $\mu(I)$ and $\varphi_{x}''(I)$ are compact.

The mapping $\gamma:I \rightarrow B_{2}$,
$$\gamma(t)=\cases{
                   \widetilde{\gamma}(t), & when $t \in [0,t_{1})
                   \cup (t_{2}, 1]$ \cr
                   \mu(\frac{t-t_{1}}{t_{2}-t_{1}}), & when $t \in
                   [t_{1}, t_{2}]$ \cr
}$$
satisfyes the following conditions: $\gamma(0)=y$, $\gamma(1)=z$,
$\gamma((0,1)) \subset \Int B_{2}$, $\gamma$ have a transversal
intersection with $\varphi_{x}''$.

Let $\{ J_{i}=\varphi_{x}'' ([t_{i,1},t_{i,2}]) , \; i \in \nn \}$ be a
system of connected segments of $\varphi_{x}''$ such that
$$\varphi_{x}''(t_{i,k}) \in \partial B_{2}, \; i=1, \ldots, n, \; k
= 1,2 \; ;$$
$$\varphi_{x}''((t_{i,1},t_{i,2})) \subset \Int B_{2}, \; i=1, \ldots, n,
\; ;$$
$$\varphi_{x}''(I) \cap \gamma(I) \subset \bigcup_{i=1}^{n} J_{i}
\;.$$

We are going to show the existance of $j \in \{ 1, \ldots, n \}$ such that
$y$ and $z$ are contained in the different arcwise connected components of
$B_{2} \setminus J_{j}$.

Assume it is not a case and the points $y$ and $z$ are situated in a same
component of $B_{2} \setminus J_{i}$ for all $i \in \{ 1, \ldots, n\}$.

Each set $J_{i}$, $i \{ \in 1, \ldots, n \}$ decomposes $B_{2}$ into two
arcwise connected components. Denote $\widetilde{B}_{2}(i)$ the component
of $B_{2} \setminus J_{i}$ containing $y$ and $z$.

Set
$$\tau_{1} = \min \{ t \; | \; \gamma(t) \in J_{1} \}$$
$$\tau_{2} = \max \{ t \; | \; \gamma(t) \in J_{1} \} \;.$$
It is obvious that $\gamma([0, \tau_{1})) \cup \gamma((\tau_{2}, 1])
\subset \widetilde{B}_{2}(1)$, $\gamma(\tau_{1}) = \varphi_{x}''(s_{1})$,
$\gamma(\tau_{2}) = \varphi_{x}''(s_{2})$ for certain $s_{1}$, $s_{2}$
such that $t_{1,1} \le s_{k} \le t_{1,2}$, $k = 1,2$. Notice that $s_{k}
\ne t_{1,k}$, $k=1,2$ since $\varphi_{x}''(s_{1})$, $\varphi_{x}''(s_{2})
\in \Int B_{2}$.

Consider compact nonintersecting sets $K_{1} = \varphi_{x}''([0, t_{1,1}]
\cup [t_{1,2}, 1])$, $K_{2} = \varphi_{x}''([s_{1}, s_{2}])$.

Denote
$$\varepsilon_{1} = \dist (K_{1}, K_{2}) = \min\limits_{x_{1} \in K_{1},
x_{2} \in K_{2}} \dist (x_{1}, x_{2}) \;.$$
It is evident that $\varepsilon_{1} > 0$.

The set $\Cl \widetilde{B}_{2}(1)$ is homeomorphic to the closed disk $D =
\{ (x,y) \in \rr^{2} \; | \; x^{2}+y^{2} \le 1 \}$. We shall take a
homeomorphism $f:\Cl \widetilde{B}_{2}(1) \rightarrow D$. Under the
mapping $f$ the set $K_{2}$ turns to an interval $\widetilde{K}_{2}
\subset \partial D$.

Mappings $f$, $f^{-1}:D \rightarrow \Cl \widetilde{B}_{2}(1)$ are
uniformly continuous. Find $\delta_{1} > 0$ such that the inequality
$\dist (x_{1}, x_{2}) < \delta_{1}$ implies $\dist (f^{-1}(x_{1}),
f^{-1}(x_{2})) < \varepsilon_{1}$ for every $x_{1}, x_{2} \in D$.
Also $\delta > 0$ could be found to realize an implication of $\dist
(f(y_{1}), f(y_{2})) < \delta_{1}$ from $\dist (y_{1}, y_{2}) < \delta$.

We find $\tau_{1}' < \tau_{1}$ and $\tau_{2}' > \tau_{2}$ complying with
the condition $\dist (\gamma(\tau_{k}'), \gamma(\tau_{k})) = \dist
(\gamma(\tau_{k}'), \varphi_{x}''(s_{k})) < \min (\delta,
\varepsilon_{1}/2)$, $k=1,2$. Then $f \circ \gamma(\tau_{1}')$,
$f \circ \gamma(\tau_{2}') \in U_{\delta_{1}}(\widetilde{K}_{2})
= \{ x \in D \; | \; \dist (x, \widetilde{K}_{2}) < \delta_{1} \}$.

We join the points $f \circ \gamma(\tau_{1}')$ and $f \circ
\gamma(\tau_{2}')$ with a continuous injective path $\widetilde{\alpha}:I
\rightarrow D$ such that $\widetilde{\alpha}(0) = f \circ
\gamma(\tau_{1}')$, $\widetilde{\alpha}(1) = f \circ \gamma(\tau_{2}')$,
$\widetilde{\alpha}(I) \subset (\Int D \cap U_{\delta_{1}}
(\widetilde{K}_{2}))$.

The mapping $\alpha = f^{-1} \circ \widetilde{\alpha}:I \rightarrow B_{2}$
posesses the following properties: $\alpha (0) = \gamma (\tau_{1}')$,
$\alpha (1) = \gamma (\tau_{2}')$, $\alpha (I) \subset \Int
\widetilde{B}_{2}(1) \cap U_{\varepsilon_{1}}(\varphi_{x}''([s_{1},
s_{2}]))$ (the latter implies $\alpha (I) \cap (\varphi_{x}''(I) \setminus
J_{1}) = \emptyset$), $\alpha (I) \cap J_{1} = \emptyset$.

Assume
$$\gamma_{1}(t) = \cases{
                   \gamma(t), & for $t \in [0, \tau_{1}')
                   \cup (\tau_{2}', 1]$ \cr
                   \alpha (\frac{t-\tau_{1}'}{\tau_{2}'-\tau_{1}'}),
                   & for $t \in [\tau_{1}', \tau_{2}']$ \cr
}$$

The mapping $\gamma_{1}:I \rightarrow B_{2}$ is continuous and the
following conditions are realized:
$$\gamma_{1}(0) = y, \; \gamma_{1}(1) = z, \; \gamma_{1}((0, 1)) \subset
\Int B_{2} \; ;$$
$$\gamma_{1}(I) \cap \varphi_{x}''(I) \subset \bigcup_{i=2}^{n} J_{i} \;
;$$
$$\gamma_{1}(I) \cap \varphi_{x}''(I) \subset \gamma(I) \cap
\varphi_{x}''(I) \;.$$

Presuppose the mapping $\gamma_{k}:I \rightarrow B_{2}$ complying with the
conditions
$$\gamma_{k}(0) = y, \; \gamma_{k}(1) = z, \; \gamma_{k}((0, 1)) \subset
\Int B_{2} \; ;$$
$$\gamma_{k}(I) \cap \varphi_{x}''(I) \subset \bigcup_{i=k+1}^{n} J_{i} \;
;$$
$$\gamma_{k}(I) \cap \varphi_{x}''(I) \subset \gamma(I) \cap
\varphi_{x}''(I) \;.$$
is already built for certain $k \in \{ 1, \ldots, n-1 \}$.

Let us build a mapping $\gamma_{k+1}:I \rightarrow B_{2}$ with analogously
properties.

Put
$$\tau_{1} = \min \{ t \; | \; \gamma_{k}(t) \in J_{k+1} \}$$
$$\tau_{2} = \max \{ t \; | \; \gamma_{k}(t) \in J_{k+1} \} \;.$$

Like beyond, we can find points $\tau_{1}' < \tau_{1}$, $\tau_{2}' >
\tau_{2}$ and a mapping $\alpha_{k}:I \rightarrow B_{2}$ such that
$\alpha_{k}(0) = \gamma (\tau_{1}')$, $\alpha_{k}(1) = \gamma (\tau_{2}')$
$\alpha_{k}(I) \subset \Int B_{2}$, $\alpha_{k}(I)
\cap \varphi_{x}''(I) = \emptyset$.

Assume
$$\gamma_{k+1}(t) = \cases{
                    \gamma_{k}(t), & for $t \in [0, \tau_{1}')
                    \cup (\tau_{2}', 1]$ \cr
                    \alpha_{k}(\frac{t-\tau_{1}'}{\tau_{2}'-\tau_{1}'}),
                    & for $t \in [\tau_{1}', \tau_{2}']$ \cr
}$$

The mapping $\gamma_{n}:I \rightarrow B_{2}$ complies with the equality
$\gamma_{n}(I) \cap \varphi_{x}''(I) = \emptyset$. But this is impossible
since $y \in \Int A_{1}'' = \Int (I \setminus A_{2}'')$, $z \in \Int
A_{2}''$ and $\partial A_{1}'' = \partial A_{2}'' = \varphi_{x}''(I)$ in
the topology induced from $I^{2}$.

Thus $j \in \{ 1, \ldots, n \}$ could be found such that the points $y$
and $z$ are contained in the distinct connected components of $B_{2}
\setminus J_{j}$.

Assume $J_{j} = \varphi_{x}''([t_{1}, t_{2}])$, $0 \le t_{1} < t_{2} <
1/2$ (the case $1/2 < t_{1} < t_{2} \le 1$ is considered analogously).

Points $y$ and $z$ decompose $\partial B_{2}$ into two connected
components, one of which is situated in $\varphi_{y}([0, \frac{1}{2}))
\cup I_{b} \cup I_{r}$, other is lying in $\varphi_{y}((\frac{1}{2}, 1])
\cup I_{t} \cup I_{r}$. Also $J_{j} \cap I_{r} = \emptyset$ and
$(\varphi_{y}([0, \frac{1}{2})) \cup I_{b}) \cap
(\varphi_{y}((\frac{1}{2}, 1]) \cup I_{t}) = \emptyset$. Consequently one
of the points $\varphi_{x}''(t_{1})$, $\varphi_{x}''(t_{2})$ is contained in
the set
$$\beta_{1} = \varphi_{y}([0, \frac{1}{2})) \cup I_{b} \;,$$
other lies in
$$\beta_{2} = \varphi_{y}((\frac{1}{2}, 1]) \cup I_{t} \;.$$

Suppose $\varphi_{x}''(t_{1}) \in \beta_{1}$, $\varphi_{x}(t_{2}) \in
\beta_{2}$. Mark $\varphi_{x}(t_{2}) \not\in I_{t}$ since $t_{2} < 1/2 <
1$.

There exists a unique $\tau_{2} \in (1/2, 1)$ such that
$\varphi_{y}(\tau_{2}) = \varphi_{x}''(t_{2})$. One of two following
possibilities could be realized: \\
1) Either $\varphi_{x}''(t_{1}) \in I_{b}$. Then $t_{1} = 0$ and the
mapping
$$\psi:I \rightarrow I^{2} \;,$$
$$\psi(t) = \cases{
                    \varphi_{x}''(t), & for $t \in [0, t_{2})$ \cr
                    \varphi_{y}(\frac{1-\tau_{2}}{1-t_{2}}(t-t_{2})+\tau_{2}),
                    & for $t \in [t_{2}, 1]$ \cr
}$$
satisfies the statement of lemma.

\noindent
2) Or $\varphi_{x}''(t_{1}) \in \varphi_{y}((0, 1/2))$. Then $t_{1} \in
(0, 1/2)$. Find $\tau_{1} \in (0, 1/2)$ to satisfy the equality
$\varphi_{y}(\tau_{1}) = \varphi_{x}''(t_{1})$. The mapping
$$\psi:I \rightarrow I^{2} \;,$$
$$\psi(t) = \cases{
                    \varphi_{y}(t), & for $t \in [0, \tau_{1})
                    \cup (\tau_{2}, 1]$ \cr
                    \varphi_{x}''(\frac{t_{2}-t_{1}}{\tau_{2}-\tau_{1}}
                    (t-\tau_{1})+t_{1}), & for $t \in
                    [\tau_{1}, \tau_{2}]$ \cr
}$$
complies with the statement of lemma.

\bigskip

Enter upon testing the correctness of the definition given above.

Let $x,y$ be arbitrary elements of $K \setminus (\{ a \} \cup \{ b \})$.

1) Let us show the definition does not depend on a choise of the mapping
$\varphi_{x}:I \rightarrow I^{2}$.

Let us take two injective continuous mappings $\varphi_{x}',
\varphi_{x}'':I \rightarrow I^{2}$ to realize the following conditions:
$\varphi_{x}'(0), \varphi_{x}''(0) \in \stackrel{\circ}{I_{b}}$,
$\varphi_{x}'(1), \varphi_{x}''(1) \in \stackrel{\circ}{I_{t}}$,
$\varphi_{x}'(\frac{1}{2}) = \varphi_{x}'' (\frac{1}{2}) = x$,
$\varphi_{x}'([0, \frac{1}{2})), \varphi_{x}''([0, \frac{1}{2})) \subset
D_{1}$,
$\varphi_{x}'((\frac{1}{2},1]), \varphi_{x}''((\frac{1}{2},1]) \subset
D_{2}$, $y$ is contained both in the same connected component of $I^{2}
\setminus \varphi_{x}'(I)$ with $I_{l}$ and in the same component of
$I^{2} \setminus \varphi_{x}''(I)$ with $I_{r}$.

In accord with proposition \ref{prop2} there exists injective continuous
mapping $\varphi_{y}:I \rightarrow I^{2}$, $\varphi_{y}(0) \in
\stackrel{\circ} {I_{b}}$, $\varphi_{y}(1) \in \stackrel{\circ}{I_{t}}$,
$\varphi_{y} (\frac{1}{2}) = y$, $\varphi_{y}([0, \frac{1}{2})) \subset
D_{1}$ $\varphi_{y}((\frac{1}{2},1]) \subset D_{2}$ such that
$\varphi_{y}(I) \cap \varphi_{x}'(I) = \emptyset$. Our construction
satisfies the conditions of lemma \ref{lemma4}. Therefore there exists a
continuous injective mapping $\psi:I \rightarrow I^{2}$ such that
$$\psi(I) \subset \varphi_{x}''(I) \cup \varphi_{y}(I) \;,$$
$$\psi(0) \in \stackrel{\circ}{I_{b}} \subset D_{1} \;,\;
\psi(1) \in \stackrel{\circ}{I_{t}} \subset D_{2} \;,$$
$$x, y \not\in \psi(I) \;.$$

Thus $\psi(I) \cap K = \emptyset$ what is impossible.

2) Let $\varphi_{x}:I \rightarrow I^{2}$ be a mapping complying with the
statement of proposition \ref{prop1}. Suppose $y$ is lying in the same
connected
component of $I^{2} \setminus \varphi_{x}(I)$ with $I_{l}$. Let us take a
mapping $\varphi_{y}:I \rightarrow I^{2}$ from proposition \ref{prop2}. An
emidiate verification shows that $x$ is contained in the same connected
component of $I^{2} \setminus \varphi_{y}(I)$ with $I_{r}$.

\bigskip

Let us take $x, y \in K$. We denote by $W(x,y)$ the set $\{ z \in K \;|\;
x \le z \le y \}$ if $x \le y$ and the set $\{ z \in K \;|\; y \le z \le x
\}$ otherwise.

\begin{lemma}
{\rm
Let $z_{1}, z_{2}$ be arbitrary elements of $K$ complying with $z_{1} \le
z_{2}$. The set $W(z_{1}, z_{2})$ is compact. Under $z_{1} \ne z_{2}$ the
set $\{ z \in K \: | \: z_{1} \le z \le z_{2}, \; z \ne z_{1},z_{2} \}$
appeares to be the nonempty opened set in the topology induced on $K$ from
$(I^{2}, \dist)$.
\label{lemma5}
}
\end{lemma}

\noindent
{\bf Proof.} For $z_{1}=z_{2}$ the set $W(z_{1}, z_{2})$ consists of a
single point, thus it is compact.

Let $z_{1} \ne z_{2}$. Under $z_{1}=a$ take the mapping $\varphi_{z_{1}}:I
\rightarrow I_{l}$, $\varphi_{z_{1}}(t) = \{ 0 \} \times \{ t \}$. Under
$z_{2}=b$ take the mapping $\varphi_{z_{2}}:I \rightarrow I_{r}$,
$\varphi_{z_{2}}(t) = \{ 1 \} \times \{ t \}$. Otherwise find a mapping
$\varphi_{z_{1}}:I \rightarrow I^{2}$ meeting the statement of proposition
\ref{prop1} and a mapping $\varphi_{z_{2}}:I \rightarrow I^{2}$ meeting the
statement of proposition \ref{prop2}.

Consider an injective continuous mapping
$$\gamma:S^{1} \rightarrow I^{2} \;,$$
$$\gamma(t)=\cases{
                    \varphi_{z_{1}}(4t) & for $t \in [0, 1/4)$ \cr
                    ((2-4t)\varphi_{z_{1}}(1)+(4t-1)\varphi_{z_{2}}(1) ,1)
                          & for $t \in [1/4, 1/2)$ \cr
                    \varphi_{z_{2}}(3-4t) & for
                          $t \in [1/2, 3/4)$ \cr
                    ((4t-3)\varphi_{z_{1}}(0)+
                          (4-4t)\varphi_{z_{2}}(0) ,0) & for
                          $t \in [3/4, 1]$. \cr
                 }$$
bounding a closed disk $D_{z_{1}, z_{2}}$, such that $D_{z_{1}, z_{2}}
\cap K = W(z_{1}, z_{2})$,
$\Int D_{z_{1}, z_{2}} \cap K = \{ z \in K \: | \: z_{1} \le z \le z_{2},
\; z \ne z_{1},z_{2} \}$.

The set $\Int D_{z_{1}, z_{2}} \cap K$ could not be empty since
$\Int D_{z_{1}, z_{2}} \cap D_{1} \ne \emptyset$,
$\Int D_{z_{1}, z_{2}} \cap D_{2} \ne \emptyset$ with $\Int D_{z_{1},
z_{2}}$ is arcwise connected.

As $K$ is closed and $D_{z_{1}, z_{2}}$ is compact, the set $D_{z_{1},
z_{2}} \cap K = W(z_{1}, z_{2})$ is compact. \\
Lemma is proved.

\bigskip

We shall define the function
$$\myrho: K \times K \rightarrow \rr_{+} \;,\; \rho(x,y) = \diam W(x,y) \;.$$

\begin{prop}
{\em
$\myrho: K \times K \rightarrow \rr_{+}$ is a distance function on $K$.
\label{prop3}
}
\end{prop}

\noindent
{\bf Proof.} $\myrho$ is a nonnegative function by definition.

Let us verify the fulfillment for $\myrho$ of the distance function
properties.

1) Obviously $\myrho (z, z)=0$ for any $z \in K$. Let $\myrho (z_{1},
z_{2})=0$ for certain $z_{1}$, $z_{2} \in K$. Since $\myrho (z_{1}, z_{2})
= \diam W(z_{1}, z_{2}) \ge \dist (z_{1}, z_{2})$, $z_{1}=z_{2}$.

2) $\myrho (z_{1}, z_{2}) = \diam W(z_{1}, z_{2}) = \myrho (z_{2}, z_{1})$
for any $z_{1}$, $z_{2} \in K$.

3) Fix $z_{1}$, $z_{2}$, $z_{3} \in K$. Let us show the correctness of
$$\myrho (z_{1}, z_{3}) \le \myrho (z_{1}, z_{2}) + \myrho (z_{2},
z_{3}) \;.
\eqno (*)$$

Assume $z_{2} \not\in W(z_{1}, z_{3})$. In this case either $W(z_{1},
z_{3}) \subset W(z_{1}, z_{2})$ or $W(z_{1}, z_{3}) \subset W(z_{2},
z_{3})$. The inequality (\mystar) is fulfilled as $B \subset A$ implies
$\diam B \ge \diam A$ for all limited $A, B \subset \rr^{2}$.

Let $z_{2} \in W(z_{1}, z_{3})$. Since the metric $\dist:I^{2} \times
I^{2} \rightarrow \rr_{+}$ is continuous and the set $W(z_{1}, z_{3})
\times W(z_{1}, z_{3}) \subset I^{2} \times I^{2}$ is compact, $(y_{1},
y_{2}) \in W(z_{1}, z_{3}) \times W(z_{1}, z_{3})$ could be found to
realize an equality
$$\dist (y_{1}, y_{2}) = \max\limits_{x_{1}, x_{2} \in W(z_{1}, z_{3})}
\dist (x_{1}, x_{2}) = \myrho (z_{1}, z_{3}) \;.$$

Under $y_{1}$, $y_{2} \in W(z_{1}, z_{2})$ or $y_{1}$, $y_{2} \in W(z_{2},
z_{3})$ either an inequality $\myrho (z_{1}, z_{2}) \ge \dist (y_{1},
y_{2}) = \myrho (z_{1}, z_{3})$ or $\dist (y_{1}, y_{2}) = \myrho (z_{1},
z_{3}) \le \myrho (z_{2}, z_{3})$ is valid and the relation (\mystar)
is satisfyed.

In the case $y_{1} \in W(z_{1}, z_{2})$, $y_{2} \in W(z_{2}, z_{3})$ we
have $\myrho (z_{1}, z_{3}) = \dist (y_{1}, y_{2}) \le \dist (y_{1},
z_{2}) + \dist (z_{2}, y_{2}) \le \myrho (z_{1}, z_{2}) + \myrho (z_{2},
z_{3})$.

The case $y_{1} \in W(z_{2}, z_{3})$, $y_{2} \in W(z_{1}, z_{2})$ can be
reduced to the preceding since $\dist$ is symmetric. \\
Proposition is proved.

\begin{defn}
{\em
Call a sequence $\{ z_{k} \in K \}_{k \in \nn}$ {\it nondecreasing}
(respectively {\it nonincreasing}) in the case $z_{k} \le z_{k+1}$ for any
$k \in \nn$ ($z_{k} \ge z_{k+1}$ for any $k \in \nn$).

Call a sequense $\{ z_{k} \in K \}_{k \in \nn}$ {\it monotonic} if it is
nondecreasing or nonincreasing.
}
\end{defn}

\begin{defn}
{\em
Call $z \in K$ an {\it upper} ({\it lower}) {\it bound} of a sequence $\{
z_{k} \in K \}_{k \in \nn}$ if $z_{k} \le z$ for any $k \in \nn$ ($z_{k}
\ge z$ for any $k \in \nn$).

Call $z \in K$ a {\it least upper} ({\it greatest lower}) {\it bound} of a
sequense $\{ z_{k} \in K \}_{k \in \nn}$ in case $z$ appears an upper
(lower) bound of $\{ z_{k} \in K \}_{k \in \nn}$ and any $z' \le z$, $z'
\ne z$ ($z' \ge z$, $z' \ne z$) is not the upper (lower) bound of this
sequence.
}
\end{defn}

\begin{lemma}
{\em
Let $\{ z_{k} \in K \}_{k \in \nn}$ be a nonincreasing (nondecreasing)
sequence. At that time it converges in $(I^{2}, \dist)$ to a certain $z
\in K$ being a greatest lower (least upper) bound of $\{ z_{k} \}$.
\label{lemma6}
}
\end{lemma}

\noindent
{\bf Proof.} Let $\{ z_{n} \in K \}$ be a nondecreasing. Since $K$ is
compact a limit point $z \in K$ could be found for this sequence. Select a
subsequence $\{ z_{n_{s}} \}$ converging to $z$.

Let us verify the inequality $z \ge z_{n}$ is valid for any $n \in \nn$.
Presuppose it is not the case. At that time $n_{0}$ could be found such
that $z \le z_{n_{0}}$, $z \ne z_{n_{0}}$, hence $z \le z_{n}$, $z \ne
z_{n}$ for all $n \ge n_{0}$. In accord with lemma \ref{lemma5} we can
find $y \in K$ meeting $z \le y \le z_{n}$, $y \ne z$, $z_{n_{0}}$.

Take a mapping $\varphi_{y}:I \rightarrow I^{2}$ complying with the
statement of proposition \ref{prop1}. It divides $I^{2}$ into two closed
disks $A_{1} \supset W(a, y)$, $A_{2} \supset W(y, b)$, $A_{1} \cap A_{2} =
\varphi_{y}(I)$. Note $z$ is contained in $\Int A_{1}$ and also $z_{n} \in
A_{2}$ implies $z_{n} \not\in \Int A_{1}$ for all $n \ge n_{0}$. So the
set $\Int A_{1}$ is an opened neighbourhood of $z$ which includes a finite
number of elements from $\{ z_{n} \}$, contradictory to a fact $z$ is the
limit point of $\{ z_{n} \}$.

Show $z$ is the least upper bound of $\{ z_{n} \}$. Assume it is not the
case and there exists $z' \ne z$ such that $z' \le z$ and $z_{n} \le z'$
for all $n \in \nn$. Take a mapping $\varphi_{z'}:I \rightarrow I^{2}$
meeting the statement of proposition \ref{prop1}. It divides $I^{2}$ into
disks $A_{1}'$, $A_{2}'$ where $z \in \Int A_{2}'$ and $y_{n} \in A_{1}'$
for all $n \in \nn$. This contradicts to the fact $z$ is the limit point of
$\{ z_{n} \}$, since $y_{n} \not\in \Int A_{2}'$.

Presuppose there exists $z' \ne z$ being another limit point for $\{ z_{n}
\}$. Iterating argument given above deduce point $z'$ is the least upper
bound of $\{ z_{n} \}$, thus $z' = z$ and $z$ is a limit of $\{ z_{n} \}$.
\\ Lemma is proved.

\begin{corr}
{\em
There exists an unumbiguously defined least upper (greatest lower) bound
of nonincreasing (nondecreasing) sequence $\{ z_{k} \in K \}_{k \in \nn}$.
}
\end{corr}

\begin{lemma}
{\em
A monotonic subsequence could be selected from any sequence $\{ z_{k} \in
K \}_{k \in \nn}$.
}\label{lemma7}
\end{lemma}

\noindent
{\bf Proof.} Since $K$ is a compact every sequence $\{ z_{n} \in K \}$ have
a convergent subsequence. Thus it will be enough to prove lemma for
convergent sequences.

Fix a sequence $\{ z_{n} \in K \}$ converging to a certain $z \in K$. As
$K = W(a, z) \cup W(z, b)$, either $W(a, z)$ or $W(z, b)$ contains an
infinit number of elements from the sequence. Suppose this is $W(a, z)$.
Without loss of generality we can consider $z_{n}$ to be contained in
$W(a, z)$, i. e. $z_{n} \le z$ for any $n \in \nn$.

If a stationary subsequence could be selected from $\{ z_{n} \}$ the
subsequence will be monotonic. Suppose this is not the case. Therefore a
subsequence $\{ z_{n_{k}} \}$ could be chosen with mutually distinct
elements. We can also assume  $z_{n} \ne z$ for any $n \in \nn$.

Mark off in our case the inequalities $z_{n_{1}} \le z_{n_{2}}$,
$z_{n_{1}} \ge z_{n_{2}}$ could not be valid simultaniously for $n_{1} \ne
n_{2}$ since $z_{n_{1}} \ne z_{n_{2}}$.

Let $n_{0}$ be a value such that $z_{n} \le z_{n_{0}}$ for all $n > n_{0}$.
Take a mapping $\varphi_{z_{n_{0}}}:I \rightarrow I^{2}$ satisfying the
statement of proposition \ref{prop1}. It divides $I^{2}$ into two closed
disks $A_{1}$ and $A_{2}$ complying with the relations $A_{1} \cap A_{2} =
\varphi_{z_{n_{0}}}(I)$, $A_{1} \cap K = W(a, z_{n_{0}})$, $A_{2} \cap K =
W(z_{n_{0}}, b)$. On one hand $z \in \{ x \in K \; | \; z_{n_{0}} < x \le
b, \; x \ne z_{n_{0}} \} \subset \Int A_{2}$ (in the topology induced from
$I^{2}$), since $z_{n_{0}} \in W(a, z)$ and $z_{n_{0}} \ne z$. On the
other hand $z_{m} \in A_{1}$ for any $m > n_{0}$, thus $z_{m} \not\in \Int
A_{2}$. Hence $\Int A_{2}$ appears to be the opened neighborhood of $z$
containing the finite number of elements from $\{ z_{n} \}$. This is
impossible since $z$ is the limit point of $\{ z_{n} \}$.

Consequently for any $n \in \nn$ there exists $m > n$ such that $z_{m} \ge
z_{n}$ and a nondecreasing subsequence of $\{ z_{n} \}_{n \in \nn}$ could
be found.

\begin{prop}
{\em
Both the restriction of metric $\dist$ onto $K$ and the distance function
$\myrho$ generate the same topology on $K$.
}\label{prop4}
\end{prop}

First we will prove one simple statement.

\begin{lemma}
{\em
Let $\{ x_{n} \}$, $\{ y_{n} \}$ be nondecreasing (nonincreasing)
sequences meeting the following condition: for any $n \in \nn$ there exist
$m_{1}$, $m_{2} \in \nn$ such that $y_{m_{1}} \ge x_{n}$, $x_{m_{2}} \ge
y_{n}$ ( $y_{m_{1}} \le x_{n}$, $x_{m_{2}} \le y_{n}$).

Let $x$ be a least upper (greatest lower) bound of $\{ x_{n} \}$. Then
$x$ is the least upper (greatest lower) bound for $\{ y_{n} \}$.
}\label{lemma8}
\end{lemma}

\noindent
{\bf Proof of lemma \ref{lemma8}.} Consider the case sequences $\{ x_{n}
\}$, $\{ y_{n} \}$ are nondecreasing. Let $z_{x}$, $z_{y}$ be the least
upper bounds respectively of $\{ x_{n} \}$ and $\{ y_{n} \}$. In accord
with conditions of lemma $z_{y}$ is the upper bound for $x_{n}$ since
$y_{n} \le z_{y}$ for all $n \in \nn$ and $m_{1}$ could be found for any
$n$ to realize the inequality $x_{n} \le y_{m_{1}}$.

Therefore $z_{y} \ge z_{x}$ as $z_{x}$ is the least upper bound of
$x_{n}$. Analogously we receive $z_{x} \ge z_{y}$, thus $z_{x} = z_{y}$. \\
Q. E. D.

\bigskip

\noindent
{\bf Proof of proposition \ref{prop4}.} For any $z_{0} \in K$ both the
collection
$$W_{\varepsilon}(z_{0}) = \{ z \in K \; | \; \myrho (z, z_{0}) <
\varepsilon \} , \; \varepsilon > 0$$
is a determining system of neighborhoods in the topology generated by
$\myrho$ and the collection
$$B_{\varepsilon}(z_{0}) = \{ z \in K \; | \; \dist (z, z_{0}) <
\varepsilon \} , \; \varepsilon > 0$$
is a determining system of neighborhoods in the topology induced on $K$
from $(I^{2}, \dist)$.

1) Verify the validity of following statement: for every $\varepsilon > 0$
there exists $\delta > 0$ such that $W_{\delta}(z_{0}) \subset
B_{\varepsilon}(z_{0})$. Really, $\myrho (z, z_{0}) = \diam W(z, z_{0})
\ge \dist (z, z_{0})$ for any $z \in K$. Therefore $\myrho (z, z_{0}) <
\varepsilon$ implies $\dist (z, z_{0}) < \varepsilon$ and we can take
$\delta = \varepsilon$.

Thus the topology on $K$ induced by $\myrho$ is stronger than topology
generated by $\dist$.

2) Inversly, for every $\varepsilon > 0$ there exists $\delta > 0$ such
that $B_{\delta}(z_{0}) \subset W_{\varepsilon}(z_{0})$.

Suppose it is not the case and $\varepsilon_{0} > 0$ could be found to
meet the condition $B_{\delta}(z_{0}) \setminus W_{\varepsilon_{0}}(z_{0})
\ne \emptyset$ for any $\delta >0$. Assume $\delta_{n} = 1/2^{n}$, $n \in
\nn$. For every $n \in \nn$ take a point $y_{n} \in B_{\delta_{n}}(z_{0})
\setminus W_{\varepsilon_{0}}(z_{0})$. We receive a sequence $\{ y_{n} \}$
such that $\dist (y_{n}, z_{0}) \rightarrow 0$ under $n \rightarrow
\infty$, however $\myrho (y_{n}, z_{0}) \ge \varepsilon_{0}$ for any
$n \in \nn$.

In accord with lemma \ref{lemma7} $\{ y_{n} \}$ contains a monotonic
subsequence. Without loss of generality we can regard $\{ y_{n} \}$ to be
nondecreasing.

There exist $x_{n} \in W(y_{n}, z_{0}) = \{ x \in K \; | \; y_{n} \le x
\le z_{0} \}$ complying with the inequality $\dist (x_{n}, z_{0}) \ge
\varepsilon_{0}/2$ since $\myrho (y_{n}, z_{0}) = \diam W(y_{n}, z_{0}) \ge
\varepsilon_{0}$ for any $n \in \nn$. According to lemma \ref{lemma6}
$z_{0}$ is the least upper bound of $\{ y_{n} \}$. Hence for an arbitrary
$n \in \nn$ there exists $m$ such that $y_{m} \ge x_{n}$.

With respect to lemma \ref{lemma7} the sequence $\{ x_{n} \}$ have a
monotonic subsequence. Assume $\{ x_{n} \}$ is monotonic. Then it is
obviously undecreasing.

The sequences $\{ x_{n} \}$, $\{ y_{n} \}$ satisfies demands of lemma
\ref{lemma8}. Therefore $z_{0}$ is the least upper bound of $\{ x_{n} \}$
and according to lemma \ref{lemma6} $\{ x_{n} \}$ must converge to $z_{0}$
in the topology induced from $(I^{2}, \dist)$.

The contradiction obtained shows us the validity of 2).

Consequently the topology on $K$ induced from $(I^{2}, \dist)$ is stronger
than the topology generated by $\myrho$. Together with 1) it gives us an
equevalence of topologies on $K$ induced by $\dist$ and $\myrho$
respectively.

\begin{lemma}
{\em
For any $z_{0} \in K \setminus \{ b \}$ and $\varepsilon > 0$ there exists
$z_{1} \in W_{\varepsilon}(z_{0}) \setminus \{ z_{0} \} = \{ z \in K \; |
\; 0 < \myrho (z, z_{0}) < \varepsilon \}$ complying with an equality
$z_{0} \le z_{1}$.

For any $z_{0} \in K \setminus \{ a \}$ and $\varepsilon > 0$ there exists
$z_{2} \in W_{\varepsilon}(z_{0}) \setminus \{ z_{0} \}$ such that
$z_{2} \le z_{0}$.
}\label{lemma9}
\end{lemma}

\noindent
{\bf Proof.} Let $z_{0}$ be an arbitrary element of $K \setminus (\{ a \}
\cup \{ b \})$.

Presuppose $\varepsilon > 0$ could be found to meet the condition $\myrho
(z, z_{0}) \ge \varepsilon$ for every $z \ge z_{0}$, $z \ne z_{0}$.
Proposition \ref{prop4} maintains the existance of an $\varepsilon_{0} >
0$ realizing the relation $\dist (z, z_{0}) \ge 2 \varepsilon_{0}$ for any
$z \ge z_{0}$, $z \ne z_{0}$.

Take a mapping $\varphi_{z_{0}}:I \rightarrow I^{2}$ meeting the statement
of proposition~\ref{prop1}. It divides $I^{2}$ into two closed disks
$A_{1}$, $A_{2}$ such that
$$A_{1} \cup A_{2} = I^{2}, \; A_{1} \cap A_{2} = \varphi_{z_{0}}(I) \;,$$
$$A_{1} \cap K = \{ z \in K \; | \; z \le z_{0} \} \;,$$
$$A_{2} \cap K = \{ z \in K \; | \; z \ge z_{0} \} \;.$$

We are going to verify the existance of an injective continuous mapping
$\alpha:I \rightarrow I^{2}$ satisfying the following conditions:
$$\alpha(0) = \varphi_{z_{0}}(t_{1}), \; \alpha(1) =
\varphi_{z_{0}}(t_{2})$$
for certain $t_{1} < 1/2 < t_{2}$;
$$\alpha(I) \subset B_{\varepsilon_{0}}(z_{0}), \; \alpha((0, 1)) \subset
\Int A_{2}.$$

There exists a homeomorphism $f:A_{2} \rightarrow I^{2}$ complying with
$f \circ \varphi_{z_{0}}(I) = I_{l}$, $f(z_{0}) \in
\stackrel{\circ}{I_{l}}$.

The continuous mapping $f^{-1}:I^{2} \rightarrow A_{2}$ is well defined.
Hence $\delta_{0} > 0$ could be found such that $\dist (x, y) < 2
\delta_{0}$ have as a consequence
$$
\dist (f^{-1}(x), f^{-1}(y)) < \varepsilon_{0}
$$
for any $x$, $y \in I^{2}$.

Assume $\delta = \min (\delta_{0}, \dist (f(z_{0}), I_{t} \cup I_{b}))$.
Denote by $\widetilde{\alpha}:I \rightarrow I^{2}$ an arc of parametrized
sircle $\partial B_{\delta}(f(z_{0}))$ contained in $I^{2}$ with endpoints
in $I_{l}$. At that time the mapping $\alpha = f^{-1} \circ
\widetilde{\alpha}:I \rightarrow A_{2}$ will meet conditions desired.

Mark the sets $\alpha(I)$ and $K$ have a trivial intersection. Denote by
$\gamma:I \rightarrow I^{2}$ the injective continuous mapping
$$\gamma(t) = \cases{
                \varphi_{z_{2}}(t) & for $t \in [0, t_{1})
                \cup (t_{2}, 1]$ \cr
                \alpha(\frac{t-t_{1}}{t_{2}-t_{1}}) & for
                $t \in [t_{1}, t_{2}]$ \cr
} \;.$$

Thus the following incompatible correlations must be valid:
$$\gamma(I) \subset D_{1} \cup D_{2}, \; \gamma(I) \cap D_{1} \ne \emptyset,
\; \gamma(I) \cap D_{2} \ne \emptyset,$$
$$\gamma(I) \cap (\partial D_{1} \cap \partial D_{2}) = \gamma(I) \cap K =
\emptyset.$$

The contradiction obtained shows the correctness of the first statement of
lemma for any $z_{0} \in K \setminus (\{ a \} \cup \{ b \})$.

The argument mentioned can be extended to the case $z_{0} = a$ under
$$\varphi_{a}:I \rightarrow I^{2}, \; \varphi_{a}(t) = \{ 0 \} \times \{ t
\}$$
$$A_{2} = I^{2}, \; f = id:A_{2} \rightarrow I^{2} \;.$$

The second statement of lemma is proved analogously.

\bigskip

Let us fix points $z_{1}$, $z_{2} \in K$, $z_{1} \ne z_{2}$, $z_{1} \le
z_{2}$.

Consider functions
$$f_{z_{s}}:K \rightarrow \rr \;, s=1,2 \;,$$
$$f_{z_{s}}(z) = \cases{
                   \myrho (z, z_{s}), & for $z \ge z_{s}$ \cr
                   - \myrho (z, z_{s}), & for $z \le z_{s}$ \cr
} \;;$$

$$f_{z_{1}, z_{2}}(z):K \rightarrow \rr \;,$$
$$f_{z_{1}, z_{2}}(z) = f_{z_{1}}(z) + f_{z_{2}}(z) \;.$$

The definition of $\myrho$ with regard for proposition \ref{prop4} have as
a consequence the continuity of $f_{z_{1}, z_{2}}$. An immediate
verification shows us $f_{z_{1}, z_{2}}$ is nondecreasing on $K$.

\begin{prop}
{\em
Let $f_{z_{1}, z_{2}}(z')=t'$, $f_{z_{1}, z_{2}}(z'')=t''$ for certain
$z'$, $z'' \in K$, $z' \le z''$. For an arbitrary $t \in [t', t'']$ there
exists $z \in W(z', z'')$ being a prototype of $t$ under the mapping
$f_{z_{1}, z_{2}}(z)=t$.
}\label{prop5}
\end{prop}

\noindent
{\bf Proof.} First we are going to show
$$f_{z_{1}, z_{2}}(K) = [- \myrho(a, z_{1}) - \myrho(a, z_{2}), \;
\myrho(z_{1}, b) + \myrho(z_{2}, b)] \;.$$
Obviously
$$f_{z_{1}, z_{2}}(a) = - \myrho(a, z_{1}) - \myrho(a, z_{2}) \;,$$
$$f_{z_{1}, z_{2}}(a) = \myrho(z_{1}, b) + \myrho(z_{2}, b) \;.$$
Hence
$$f_{z_{1}, z_{2}}(K) \subset [- \myrho(a, z_{1}) - \myrho(a, z_{2}), \;
\myrho(z_{1}, b) + \myrho(z_{2}, b)] = J \;.$$

Mark $f_{z_{1}, z_{2}}(K)$ is a closed subset of the interval $J$ since
$K$ is compact and  $f_{z_{1}, z_{2}}$ is continuous.

Let $\tau \in J \setminus f_{z_{1}, z_{2}}(K)$. An interval $(t_{1},
t_{2}) \in J \setminus f_{z_{1}, z_{2}}(K)$ could be found to comply with
ratio $\tau \in (t_{1}, t_{2})$; $t_{1}, t_{2} \in f_{z_{1}, z_{2}}(K)$.

The set $K$ decomposes into two nonintersecting subsets
$$
\begin{array}{l}
K_{1} = f^{-1}_{z_{1}, z_{2}}([ - \myrho(a, z_{1}) - \myrho(a, z_{2}),
t_{1}]) \,, \vphantom{\sum\limits_{1}}\\
K_{2} = f^{-1}_{z_{1}, z_{2}}([t_{2}, \myrho(z_{1}, b) +
\myrho(z_{2}, b)]) \,,
\end{array}
$$
being closed in $K$. Sets $K_{1}$, $K_{2}$ are compact
since $K$ is compact.

Thus
$$\varepsilon = \myrho (K_{1}, K_{2}) = \min\limits_{x_{1} \in K_{1}, \;
x_{2} \in K_{2}} \myrho (x_{1}, x_{2}) > 0 \;.$$

As $f_{z_{1}, z_{2}}$ is monotonic, $z^{(1)} \le z^{(2)}$ for any $z^{(1)}
\in K_{1}$, $z^{(2)} \in K_{2}$. In accord with lemma \ref{lemma9} for an
arbitrary $x \in K_{1}$ we can found $y = y(x) \in B_{\varepsilon/3}(x)
\cap \{ z \in K \; | \; z \ge x, \; z \ne x \}$.

Verify $y(x) \not\in K_{2}$. Really, for every $z \in K_{2}$
$$ \myrho (y, z) + \myrho (y, x) \ge \myrho (z, x) \ge \varepsilon \;,$$
$$\myrho (y, z) \ge \varepsilon - \myrho (y, x) \ge
\varepsilon - \frac{\varepsilon}{3} = \frac{2}{3} \varepsilon$$
hence
$$\myrho (y, K_{2}) = \min\limits_{z \in K_{2}} \myrho (y, z) \ge
\frac{2}{3} \varepsilon > 0 \;.$$

Consequently $y(x) \in K_{1}$. And what is more an opened set $V(x) = \{ z
\in K \; | \; z \le y(x), \; z \ne y(x) \}$ is contained in $K_{1}$.

Since $x \in V(x)$ for any $x \in K_{1}$, the system $\{ V(x) \}_{x \in
K_{1}}$ form an opened covering of $K_{1}$. Select a finite subcovering
$\{ V(x_{i}) \}_{i=1}^{n}$.

Note $y(x_{i}) \le y(x_{j})$ implies $V(x_{i}) \subset V(x_{j})$ for any
$i$, $j \in \{ 1, \ldots, n \}$, $i \ne j$. Take $j$ so that $V(x_{i})
\subset V(x_{j})$ for all $i \ne j$. Then $K_{1} \subset V(x_{j}) = \{ z
\in K \; | \; z \le y(x_{j}), \; z \ne y(x_{j}) \}$ and $y(x_{j}) \not\in
K_{1}$.

We have a contradiction with $J \setminus f_{z_{1}, z_{2}}(K) \ne
\emptyset$. Therefore $f_{z_{1}, z_{2}} = J$.

Consider certain $z'$, $z'' \in K$, $z' \le z''$. Let now $f(z')=t'$,
$f(z'')=t''$, $t \in [t', t'']$. For $t=t'$ ($t=t''$) assume $z=z'$
($z=z''$ respectively).

Let $t' < t < t''$. Since $f_{z_{1}, z_{2}}$ is monotonic and $a \le z \le
b$ for any $z \in K$, then $[t', t''] \subset J$ and $t \in J$. Therefore
$z \in K$ could be found such that $t=f_{z_{1}, z_{2}}(z)$. Recalling
again the monotonicity of $f_{z_{1}, z_{2}}$ we obtain $z \in W(z_{1},
z_{2})$.

\begin{lemma}
{\em {\bf (the main).}
The set $K \subset I^{2}$ is a homeomorphic transform of an interval.
\label{lemma10}
}
\end{lemma}

\noindent
{\bf Proof.} Consider a denumerable everywhere dense subset
$$Q = \{ \frac{m}{2^{n}} \; | \; m = 1, \ldots, 2^{n}-1, 2^{n} \}$$
of the interval $I = [0, 1]$.

Construct a denumerable everywhere dense subset $Q_{K} = \{ z_{q} \}_{q
\in Q}$ of $K$. Put $z_{0} = a$, $z_{1} = b$.

Let $z_{q}$ be allready defined for all $q = m / 2^{k}$; $m = 1, \ldots ,
2^{k}-1$; $k = 1, \ldots , n-1$. In accord with proposition \ref{prop5} we
can select
$$z_{\frac{2m+1}{2^{n}}} \in \{ z \in K \; | \; \myrho (z,
z_{\frac{m}{2^{n-1}}}) = \myrho (z, z_{\frac{m+1}{2^{n-1}}}) \} =$$
$$= \{ z \in
K \; | \; f_{z_{\frac{m}{2^{n-1}}}, z_{\frac{m+1}{2^{n-1}}}}(z) = 0 \}
\subset W(z_{\frac{m}{2^{n-1}}}, z_{\frac{m+1}{2^{n-1}}}) \;,$$
for every $m \in \{ 1, \ldots , 2^{n-1}-1$ \}.

Denote
$$W_{k, n} = W(z_{\frac{k}{2^{n}}}, z_{\frac{k+1}{2^{n}}}), \; k = 1,
\ldots, 2^{n}-1, \; n \in \nn \;;$$
$$a_{n} = \max\limits_{k=1, \ldots, 2^{n} - 1} \diam (W_{k, n}) \;.$$

For any $n \in \nn$ the inequality $a_{n+1} \le a_{n}$ is valid since
$$W(z_{\frac{k}{2^{n}}}, z_{\frac{k+1}{2^{n}}}) \subset
W(z_{\frac{m}{2^{n-1}}}, z_{\frac{m+1}{2^{n-1}}}) \;,$$
where $m = [\frac{k}{2}]$, hence $\diam (W_{k, n}) \le \diam (W_{m,
n-1})$.

Show $a_{n} \rightarrow 0$ under $n \rightarrow \infty$. Suppose it is not
the case and there exist $\varepsilon > 0$ and sequence $\{ n_{i} \in \nn
\}$ ($n_{i} \rightarrow \infty$ for $i \rightarrow \infty$) such that
$a_{n_{i}} \ge \varepsilon$. The latter implies $a_{n} \ge \varepsilon$
for any $n \in \nn$.

Find a sequence $\{ W_{k_{n}, n} \; | \; \diam W_{k_{n}, n} \ge
\varepsilon \}_{n \in \nn}$. Select a sequence $\{ y_{n} \}$ with mutually
disjoint elements contained in $\{ z_{\frac{k_{n}}{n}} \}_{n \in \nn} \cup
\{ z_{\frac{k_{n}+1}{n}} \}_{n \in \nn}$ to comply with $\myrho (y_{n},
y_{m}) \ge \varepsilon$ for every $m$, $n \in \nn$, $n \ne m$. Assume
$y_{n} = z_{\frac{s_{n}}{n}}$ where
$$s_{n} = \cases{
                k_{n}, & if $k_{n}$ is odd \cr
                k_{n} + 1, & if $k_{n}$ is even \cr
} \;.$$
We have $y_{m} \ne y_{n}$ under $m > n$, since $y_{m} =
z_{\frac{s}{2^{m}}}$ and $s$ is odd (thus there is no $l \in \nn$ such that
$\frac{s_{m}}{2^{m}} = \frac{l}{2^{n}}$).

Verify the validity of an inequality $\myrho (y_{n}, y_{m})
\ge \varepsilon$ for any $n \in \nn$, $m = 1, \ldots, n-1$. Note $\myrho
(z_{\frac{s_{n}-1}{2^{n}}}, y_{n}) = \myrho (z_{\frac{s_{n}+1}{2^{n}}},
y_{n}) \ge \varepsilon$. Let $z_{\frac{s_{m}}{2^{m}}} = y_{m} \le y_{n}$.
Then $y_{m} \le z_{\frac{s_{n}-1}{2^{n}}} \le y_{n}$ and $\myrho (y_{m},
y_{n}) \ge \myrho (\frac{z_{s_{n}+1}}{2^{n}}, y_{n}) \ge \varepsilon$.
The case $y_{m} \ge y_{n}$is considered analogously.

Thus $y_{n}$ is a denumerable subset of $K$ which has no limit points. The
last is contradictory to the compactness of $K$.

So $a_{n} \rightarrow 0$ under $n \rightarrow \infty$.

(1) The set $Q_{K}$ is dense in $K$ since a set
$$\{ z_{\frac{k}{2^{n}}} \; | \; k=1, \ldots, 2^{n}-1, 2^{n} \}$$
forms an $\varepsilon$-network in $K$ as $a_{n} < \varepsilon$.

(2) The mapping $\varphi:Q_{K} \rightarrow \rr_{+}$, $\varphi(z_{r}) = r$
is monotonic. Really $z_{\frac{m_{1}}{2^{n}}} \le z_{\frac{m_{2}}{2^{n}}}$
under $m_{1} \le m_{2}$ for any $n \in \nn$. Let
$$z_{\frac{m_{1}}{2^{n_{1}}}} \le z_{\frac{m_{2}}{2^{n_{2}}}} \;.$$
For $n_{1} \le n_{2}$
$$\frac{m_{1}}{2^{n_{1}}} = \frac{2^{n_{2}-n_{1}}m_{1}}{2^{n_{2}}} \le
\frac{m_{2}}{2^{n_{2}}} \;.$$
For $n_{1} \ge n_{2}$
$$\frac{m_{1}}{2^{n_{1}}} \ge \frac{2^{n_{1}-n_{2}}m_{2}}{2^{n_{1}}} =
\frac{m_{2}}{2^{n_{2}}} \;.$$

The subsequent argument are based on two lemmas from [2].

\begin{lemma}
{\em
Let a certain subset $F_{t}$ of a set $X$ such that
\begin{itemize}
        \item[( .)] $F_{t} \subset F_{s}$ under $t < s$;
        \item[(b.)] $\cup \{ F_{t} \; : \; t \in D \} = X$
\end{itemize}
be defined for any element $t$ of an everywhere dense subset $D$ of
$\rr_{+}$.

Assume $f(x) = \inf \{ t \; : \; x \in F_{t} \}$ for every $x \in X$.
Then both $\{ x \; : \; f(x) < s \} = \cup \{ F_{t} \; : \; t \in D, \; t
< s \}$ and $\{ x \; : \; f(x) \le s \} = \cap \{ F_{t} \; : \; t \in D,
\; t > s \}$ for any $s \in \rr$.
}\label{lemma11}
\end{lemma}

\begin{lemma}
{\em
Let a certain opened subset $F_{t}$ of a topological space $X$ such that
\begin{itemize}
        \item[( .)] $\Cl F_{t} \subset F_{s}$ under $t < s$;
        \item[(¡.)] $\cup \{ F_{t} \; : \; t \in D \} = X$
\end{itemize}
be defined for any element $t$ of an everywhere dense subset $D$ of
$\rr_{+}$.

Then the function $f$ defined as $f(x) = \inf \{ t \; : \; x \in F_{t} \}$
is continuous.
}\label{lemma12}
\end{lemma}

Consider a denumerable everywhere dence subset
$$D = \{ \frac{k}{2^{n}} \; | \; k, n \in \nn \}$$
of $\rr_{+}$ and a collection of opened subsets
$$F_{t} = \cases{
                \{ z \in K \; | \; a \le z \le z_{t} , \; z \ne z_{t}
                \}, & for $t \in D \cap [0, 1]$ \cr
                K, & for $t \not\in [0, 1]$ \cr
} \;.$$

Obviously
$$\Cl F_{t} \subset \{ z \in K \; | \; a \le z \le z_{t} \}, \quad \mbox
{for } t \in D \cap [0, 1] \;;$$
$$\Cl F_{t} = K, \quad \mbox {for } t \not\in [0, 1] \;.$$

Consequently $\Cl F_{t} \subset F_{s}$ as $t < s$. Furthemore
$$\bigcup_{t \in D} F_{t} = K \;.$$

So the collection $\{ F_{t} \}_{t \in D}$ meets the conditions of lemma
\ref{lemma12}, thus the mapping
$$f:K \rightarrow \rr_{+} \;,$$
$$f(x) = \inf \{ t \; : \; x \in F_{t}\}$$
is continuous.

Since $F_{t} = K$ under $t > 1$, $f(K) \subset [0, 1]$ and $f$ can be
considered as a mapping $f:K \rightarrow [0, 1]$.

Let us verify $f:K \rightarrow [0, 1]$ is one-to-one mapping.

Let $t$ be an arbitrary element of $Q$. In accord with lemma \ref{lemma11}
we have $\{ x \; : \; f(x) = t \} = ( \cap \{ F_{s} \; : \; s \in D,
\; s > t \}) \setminus (\cup \{ F_{s} \; : \; s \in D, \; s < t \})$.
Mark $z_{t} \not\in F_{s}$ for every $s < t$, $s \in D$ and $z_{t} \in
F_{s}$ for any $s > t$, $s \in D$. Therefore $z_{t} \in \{ x \; : \; f(x)
= t \}$ and the mapping $f|_{Q_{K}}$ coincides with the mapping $\varphi :
Q_{K} \rightarrow Q$.

Show the mapping $f$ is injective.

Take arbitrary $z'$, $z'' \in K$; $z' \le z''$; $z' \ne z''$. In accord
with lemma \ref{lemma9} the opened sets
$$U_{1} = \{ z \in K \; | \; \myrho (z, z') < \frac{1}{2} \myrho (z', z'')
\} \cap \{ z \in K \; | \; z \ge z', \; z \ne z' \}$$
$$U_{2} = \{ z \in K \; | \; \myrho (z, z'') < \frac{1}{2} \myrho (z', z'')
\} \cap \{ z \in K \; | \; z \le z'', \; z \ne z' \}$$
are not empty. Since the set $Q_{K}$ is dence in $K$, we can select points
$z_{t_{1}} \in U_{2} \cap Q_{K}$ and $z_{t_{2}} \in U_{2} \cap Q_{K}$.
Therefore $z' \le z_{t_{1}} \le z_{t_{2}} \le z''$ and $z_{t_{1}} \ne
z_{t_{2}}$. The latter has as consequence
$$f(z') \le f(z_{t_{1}})=t_{1} < f(z_{t_{2}})=t_{2} \le f(z'') \;,$$
that is $f(z') \ne f(z'')$.

Demonstrate the mapping $f$ is surjective.

Take an arbitrary $\tau \in [0, 1]$. There exists a nonincreasing sequence
$\{ s_{n} \in Q \}$ such that
$$\tau = \lim\limits_{n \rightarrow \infty} s_{n} \;.$$

Consider a sequence
$$\{ z_{s_{n}} \}_{n \in \nn} \;.$$
It appears to be nonincreasing. In accord with lemma \ref{lemma6}
$\{ z_{s_{n}} \}$ converges to a sertain $z \in K$ being the greatest
lower bound of this sequence.

Since $z \le z_{s_{n}}$ for any $n \in \nn$, an equality
$$f(z) = \inf \{ t \; | \; z \in F_{t} \} \le f(z_{s_{n}}) $$
is valid and $f(z) \le \tau$. Presuppose $f(z) < \tau$. As $Q$ is dence in
$I$ there exists $t \in Q$ complying with $f(z) < t < \tau$. Hence both
$z \le z_{t}$, $z \ne z_{t}$ and $z_{t} \le z_{s_{n}}$ for any $n \in
\nn$, contradictory to the fact $z$ is a greatest lower bound of
$\{ z_{s_{n}} \}$.

It is well known that a continuous one-to-one mapping of a compact is a
homeomorphism. Thus $f:K \rightarrow [0, 1]$ appears to be a
homeomorphism. \\
Q. E. D.

So the sets $K^{(1)}$, $K^{(2)}$ are intervals intesecting by endpoints.
Consequently $\partial D = K^{(1)} \cup K^{(2)}$ is homeomorphic to a
sircle which is known to bound a disk.

\noindent
Theorem is proved.

\bigskip

{\it Literature}

\begin{enumerate}
        \item[1.] {\it Hurevicz W., Wallman H.} Dimension theory,
                Princeton, 1941.
        \item[2.] {\it Kelley John L.} General topology, Princeton, 1957.
\end{enumerate}
\end{document}